\documentclass[12pt]{article}
\usepackage[cp1251]{inputenc}
\usepackage[english,russian,ukrainian]{babel}
\usepackage{latexsym,amsfonts,amsthm,amssymb,amsmath}
\usepackage{graphicx,graphics,hhline}
\usepackage{euscript}
\usepackage[all]{xy}
\theoremstyle{plain}
\newtheorem{theorem}{Theorem}
\newtheorem{proposition}{Proposition}
\newtheorem{lemma}{Lemma}
\newtheorem{corollary}{Corollary}
\theoremstyle{definition}
\newtheorem{definition}{Definition}
\newtheorem{example}{Example}

\begin{document}

\Large

\noindent MSC 32F17, 52A30
\vskip 2mm
\smallskip

\noindent{\bf T.M.~Osipchuk (Т.М.~Осіпчук)}  

\noindent osipchuk@imath.kiev.ua

\vskip 2mm

\noindent Institute of Mathematics NAS of Ukraine

\vskip 4mm

\noindent{\bf ON SEMICONVEX SETS IN THE PLANE}

\vskip 8mm

\small

\noindent The present work considers the properties of classes of generally convex sets in the plane known as $1$-semiconvex and weakly $1$-semiconvex. More specifically, the examples of open and closed weakly $1$-semiconvex but non $1$-semiconvex sets with smooth boundary in the plane are constructed. It is proved that such sets consist of minimum four connected components. In addition, the example of closed, weakly $1$-semiconvex, and non $1$-semiconvex set in the plane consisting of three connected components is constructed. It is proved that such a number of components is minimal for any closed, weakly $1$-semiconvex, and non $1$-semiconvex set in the plane.
\vskip 4mm

\noindent В даній роботі вивчаються властивості класів узагальнено опуклих множин на площині, які називаються $1$-напівопуклими  та слабко $1$-напівопуклими. Зокрема, побудовано приклади відкритих і замкнених слабко $1$-напівопуклих, але не $1$-напівопуклих множин з гладкою межею на площині та доведено, що такі множини складаються мінімум з чотирьох компонент зв'язності. Також побудовано приклад замкненої слабко $1$-напівопуклої, але не $1$-напівопуклої множини на площині, яка складається з трьох компонент зв'язності, та доведено, що така кількість компонент є мінімальною для довільної замкненої слабко $1$-напівопуклої, але не $1$-напівопуклої множини на площині. \vskip 4mm

\noindent В данной работе изучаются свойства классов обобщенно выпуклых множеств на плоскости, которые называются $1$-полувыпуклыми и слабо $1$-полувыпуклыми. В частности, построены примеры открытых и замкнутых слабо $1$-полувыпуклых, но не $1$-полувыпуклых множеств с гладкой границей на плоскости и доказано, что такие множества состоят минимум из четырёх компонент связности. Также построен пример замкнутого слабо $1$-полувыпуклого но не полувыпуклого множества на плоскости, состоящего из трёх компонент связности, и доказано, что такое количество компонент  является минимальным для любого замкнутого слабо $1$-полувыпуклого но не $1$-полувыпуклого множества на плоскости.

\vskip 2mm


\normalsize


\section{Introduction and auxiliary results}

A class of $m$ - semiconvex sets is one of the classes of generally convex sets.  A semiconvexity
notion was proposed by Yu.~Zeliskii. 

\begin{definition}\label{def1} (\cite{Zel2}) A set $E\subset\mathbb{R}^n$ is called {\it m-semiconvex with respect to a point} $x\in \mathbb{R}^n\setminus E$,  $1\le m<n$, if there exists
an m-dimensional half-plane $L$ such that $x\in L$ and $L\cap E=\emptyset$.
\end{definition}

\begin{definition}\label{def2} (\cite{Zel2}) A set $E\subset\mathbb{R}^n$ is called {\it m-semiconvex },  $1\le m<n$, if it is m-semiconvex with respect to every point $x\in \mathbb{R}^n\setminus E$.
\end{definition}

 We shall use the following standard notations. For a set $G\subset\mathbb{R}^n$ let $\overline{G}$ be its closure, $\mathrm{Int}\, G$ be its interior, and $\partial G=\overline{G}\setminus\mathrm{Int}\, G$ be its boudary.

We say that a set $A$ {\it is approximated from the outside} by a family of open sets $A_k$,
$k=1,2,\ldots$, if $\overline{A}_{k+1}$
is contained in $A_k$, and
$A=\cap_kA_k$ (see \cite{Aiz3}).

Yu.~Zeliskii also suggested to distinguish $m$-semiconvex and weakly
$m$-semiconvex sets.

\begin{definition}\label{def3} (\cite{Zel3}) An open set  $G\subset\mathbb{R}^n$ is called {\it
weakly $m$-semiconvex},  $1\le m<n$, if it is $m$-semiconvex for any point  $x\in\partial G$. A closed set $E\subset\mathbb{R}^n$ is called {\it weakly $m$-semiconvex} if it can be
approximated from the outside by a family of open weakly $m$-semiconvex sets.
\end{definition}

Thus, any weakly $m$-semiconvex set  $A$ is either open or closed. Among closed weakly $m$-semiconvex sets there also are sets with empty interior:
$$
A=\overline{A}=\overline{A}\setminus \mathrm{Int}\, A=\partial A.
$$
.

\begin{theorem}\label{lemm1} \textup{(\cite{Zel3} )} Let a set $E\subset \mathbb{R}^2$ be open, weakly $1$-semiconvex and non $1$-semiconvex. Then set $E$ is disconnected.
\end{theorem}

The formulation of the following theorem is equivalent to Theorem \ref{lemm1} but will be also in use.

\begin{theorem}\label{corol1}
Let a set $E\subset \mathbb{R}^2$ be open, connected and weekly $1$-semiconvex. Then $E$ is $1$-semiconvex.
\end{theorem}

The maximal connected subsets $A_i$, $i=1,2,\ldots$, of a nonempty set $A\subset\mathbb{R}^n$ are called {\it the components} of set $A$ (see \cite{Aiz3}). Herewith, $A=\cup_iA_i$.

In \cite{Zel3} it was made the assumption that any open, weakly $1$-semiconvex, and non $1$-semiconvex set consists of not less than three components. This proposition was proved in \cite{Dak9}.

\begin{theorem}\label{lemm11} \textup{(\cite{Dak9} )} Let a set $E\subset \mathbb{R}^2$ be open, weakly $1$-semiconvex, and non $1$-semiconvex. Then $E$ consists of not less than three components.
\end{theorem}

The present work proceeds  the research of Yu. Zelinskii  by investigating properties of
$1$-semiconvex and weakly $1$-semiconvex closed sets and open (closed) sets with smooth boundary in the plane.

 \section{Main results}

Everywhere further, the interval between points $x,y\in\mathbb{R}^n$ will be denoted as $xy$ and the distance between the points  will be denoted as $|x-y|$.

Let us provide a number of accessory propositions.

\begin{proposition}\label{pro1}  If a set $E\subset\mathbb{R}^n$ is $m$-semiconvex, then it  is weakly $m$-semiconvex.
\end{proposition}

\begin{proposition}\label{pro3}  There exist weakly $1$-semiconvex sets in the plane which are non $1$-semiconvex.
\end{proposition}

The following sets are an example of open, weakly $1$-semiconvex, and non $1$-semiconvex sets.

\begin{example}\label{example1}
Let  $E_1:=A_1A_2A_3A_4$  be an open rectangle with sides parallel to the axes. Let $E_2:=B_1B_2\ldots B_6$ be an open set as on Figure 1\,a.  Let $B_0\in B_1O$, where $O$ is the origin.  Let us consider points $B_t\in B_0B_1$, $t\in[0,1]$, defined as $B_t:=tB_1+(1-t)B_0$. Let $\gamma_t$, $t\in[0,1]$, be the ray starting at $B_t$ and passing through the point $A_1$.  Let $\widetilde{\gamma}_t$  be the ray symmetric to  $\gamma_t$ with respect to the axis $Oy$. Then $\widetilde{\gamma}_t$ starts at a point $\widetilde{B}_t$ and $C_1^t:=\gamma_t\cap\widetilde{\gamma}_t$, $t\in[0,1]$. Let the length of $B_1B_6$ is such that the line $\eta$ passing through the points $A_2$, $B_6$ intersects triangle $B_0\widetilde{B}_0C_1^0$. Then  $C_2^t:=\eta\cap\widetilde{\gamma}_t$ and $C_3^t:=\eta\cap\gamma_t$, $t\in[0,1]$.

Let $\xi$ be a straight parallel to the axis $Ox$,  intersecting rays $\gamma_t$, $\widetilde{\gamma}_t$ , $t\in[0,1]$, and not intersecting $E_1$, $E_2$. Now we construct open rectangles $E^t_3:=D^t_1D^t_2D^t_3D^t_4$, $t\in[0,1]$, tangent to $\gamma_t$, $\widetilde{\gamma}_t$ and lying under the line $\xi$.

Then, by the construction, any ray starting at a point  of the interior of the triangle $C_1^tC_2^tC_3^t$ intersects set $E^t:=E_1\cup E_2\cup E_3^t$ and any ray starting at $\partial E^t$ does not intersect $E^t$, $t\in[0,1]$. Thus, each set $E^t$, $t\in[0,1]$, is weakly $1$-semiconvex and non $1$-semiconvex.
 \begin{center}
\includegraphics[width=14 cm]{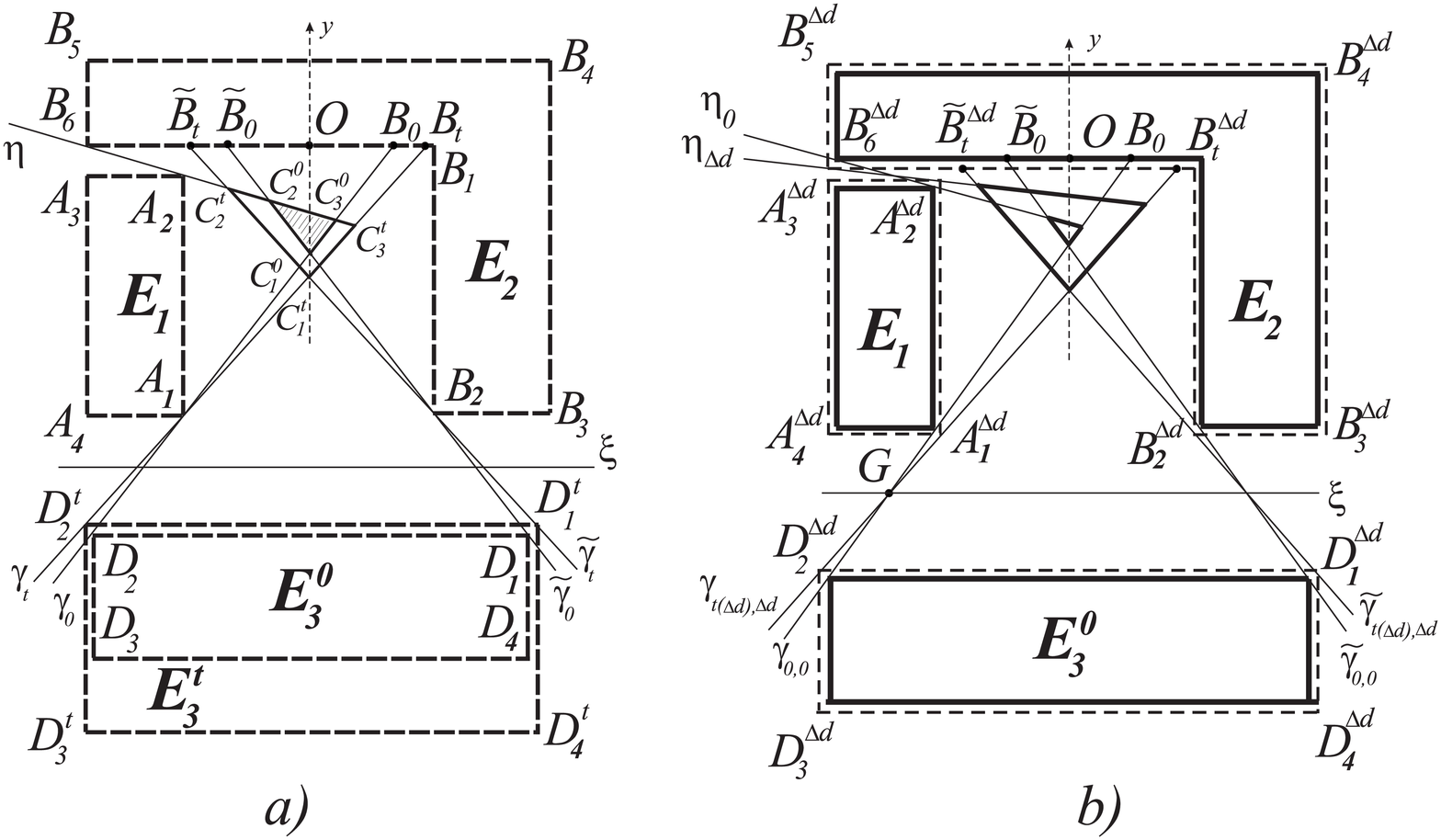}
\end{center}
\small\begin{center} Fig. 1 \end{center}
\normalsize

\end{example}

Let us construct an example of closed, weakly $1$-semiconvex, and non $1$-semiconvex set.

\begin{example}\label{example2}
Let $E^0:=E_1\cup E_2\cup E_3^0$ be the set from the previous example. The closed set $\overline{E}^0$ is non $1$-semiconvex, since any ray starting at a point of closed triangle  $C_1^0C_2^0C_3^0$ intersects $\overline{E}^0$.

Let $E_1^{\Delta d}:=A^{\Delta d}_1A^{\Delta d}_2A^{\Delta d}_3A^{\Delta d}_4\supset \overline{E}_1$, $\Delta d>0$, be the rectangle such that $|A^{\Delta d}_1A^{\Delta d}_2|={\Delta d}+|A_1A_2|$, $|A^{\Delta d}_1A^{\Delta d}_4|={\Delta d}+|A_1A_4|$ and  $E_1^{\Delta d_2}\subset E_1^{\Delta d_1}$ for any ${\Delta d_2}<\Delta d_1$, Figure 1\,b.

Let $E^{\Delta d}_2:=B^{\Delta d}_1B^{\Delta d}_2\ldots B^{\Delta d}_6\supset \overline{E}_2$ and $E_2^{\Delta d_2}\subset E_2^{\Delta d_1}$ for any ${\Delta d_2}<\Delta d_1$. Let $O^{\Delta d}\in B^{\Delta d}_1B^{\Delta d}_6\cap Oy$ and $B^{\Delta d}_0\in B^{\Delta d}_1O^{\Delta d}_2$.  Let us consider points $B^{\Delta d}_t\in B^{\Delta d}_0B^{\Delta d}_1$, $t\in[0,1]$, defined as $B^{\Delta d}_t:=tB^{\Delta d}_1+(1-t)B^{\Delta d}_0$. Let $\gamma_{t,{\Delta d}}$, $t\in[0,1]$, be the ray starting at $B^{\Delta d}_t$ and passing through the point $A^{\Delta d}_1$.

Let $\Delta d_0$ be such that $E_1^{\Delta d_0}\cap E_2^{\Delta d_0}=\emptyset$ and $G=\gamma_{1,{\Delta d_0}}\cap \gamma_{0,0}\ne \emptyset$. For a fixed $\Delta d\le\Delta d_0$ among all rays $\gamma_{t,{\Delta d}}$, $t\in[0,1]$, we choose the one $\gamma_{t(\Delta d),{\Delta d}}$  that is passing through the point $G$.

Let $\widetilde{\gamma}_{t(\Delta d),{\Delta d}}$ be the ray symmetric to  $\gamma_{t(\Delta d),{\Delta d}}$ with respect to the axis $Oy$, $0\le\Delta d\le\Delta d_0$.

Let a straight $\xi$ from the previous example be passing through the point $G$. We construct open rectangles $E^{\Delta d}_3$, $0<\Delta d\le\Delta d_0$, tangent to rays $\gamma_{t(\Delta d),{\Delta d}}$, $\widetilde{\gamma}_{t(\Delta d),{\Delta d}}$, lying under the line $\xi$, and such that $E_3^{\Delta d_2}\subset E_3^{\Delta d_1}$ for any ${\Delta d_2}<\Delta d_1$.

Now we consider the family of open weakly $1$-semiconvex sets $E^{\Delta d}:=\bigcup\limits_{j=1}^3E^{\Delta d}_j$, $0<\Delta d\le\Delta d_0$. By the constructions, $E^{\Delta d_2}\subset E^{\Delta d_1}$ for any $0<{\Delta d_2}<\Delta d_1\le\Delta d_0$. With that $\bigcap \limits_{0<\Delta d\le\Delta d_0}E^{\Delta d}=\overline{E}^0$. Thus, the closed set $\overline{E}^0$ is weakly $1$-semiconvex and non $1$-semiconvex.

\end{example}

\begin{example}\label{example3}

The following systems of open balls are examples of weakly $1$-semiconvex and non $1$-semiconvex sets with smooth boundary.

Let $B(o_1,r)$, $B(o_2,r)$ be open balls with centers $o_1$, $o_2$ placed symmetric with respect to the origin $O$ on the axis $Ox$ (see Figure 2\,a). Let $\gamma$ be the common tangent line to the balls passing through the origin.  Suppose $x=x(t)$, $y=y(t)$ is one-to-one continuous mapping of the interval $[0,1]$ onto the arc $\breve{E_{0}E_1}\subset \partial B(o_1,r)$ such that $(x(0), y(0))\equiv E_{0}\in\gamma\cap\partial B(o_1,r)$, $(x(1), y(1))\equiv E_1\in Oo_1\cap \partial B(o_1,r)$. Let $\gamma_t$, $t\in[0,1]$, be the ray starting at the point $(x(t), y(t))\in \breve{E_{0}E_1}$ and tangent to the ball $B(o_2,r)$ from the inside, i.e intersecting $Ox$. Let $\xi$ be a line parallel to the axis $Ox$ and not intersecting the balls but intersecting rays $\gamma_t$.
 To fix the point $O_1\in Oy$, let us draw the ball $B(O_1,R_1)$  tangent to the ray $\gamma_1$ and line $\xi$. Now we consider balls $B(O_1,R_t)$, $t\in[0,1]$,  with center at the point $O_1$ and  tangent to rays $\gamma_t$.  It is clear that $R_{t^1}<R_{t^2}$ for any $t^1,t^2\in[0,1]$ such that $t^1<t^2$. And let $B(O_2,R_t)$, $t\in[0,1]$, be the balls symmetric to the corresponding balls $B(O_1,R_t)$ with respect to the origin. Then each  system of four open balls $B_t:=\{B(o_i,r), B(O_i,R_t),  i=1,2\}$, $t\in(0,1]$, is weakly $1$-semiconvex and non $1$-semiconvex set. Indeed,  by the constructions, for any boundary point of $B_t$ there exists a ray starting at this point and not intersecting the set and any ray starting at a  point of the interior of rhombus $ABCD$ generated by the intersection of ray $\gamma_t$, $t\in(0,1]$, its symmetric with respect to the axis $Ox$ ray, and their  symmetric with respect to the axis $Oy$ rays intersects the set.
\begin{center}
\includegraphics[width=14 cm]{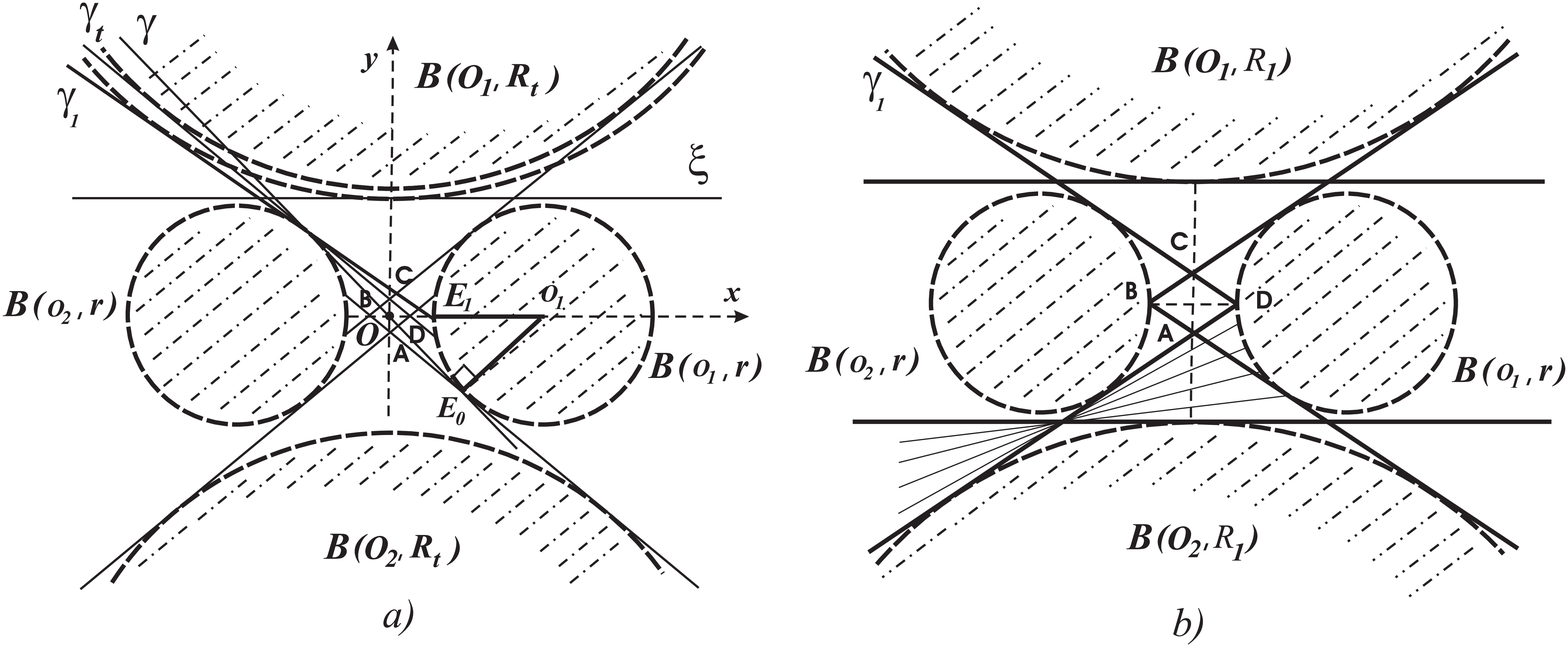}
\end{center}
\small\begin{center} Fig. 2 \end{center}

\normalsize

\end{example}

Let us construct an example of closed, weakly $1$-semiconvex, and non $1$-semiconvex set with smooth boundary.
\begin{example}\label{example4}
Let $B(o_1,r)$, $B(o_2,r)$ be open balls from the preview example and let $|o_1o_2|=2d$. Let
$\Delta r_0$ be a number such that  $0< \Delta r_0 <\dfrac{r(d-r)}{d+r}$ and let us draw a system
of concentric balls $\left\{ B(o_i,r+\Delta r),\, 0<\Delta r\le \Delta r_0,\, i=1,2\right\}$.   Let
us fix $\Delta r$ and for every ball $B(o_i,r+\Delta r)$ construct rays $\gamma_{t,\Delta r}$  as
previews. Since $\Delta r_0 <\dfrac{r(d-r)}{d+r}$, rays $\gamma_{1,\Delta r_0}$ and $\gamma_{0,0}$
are intersected at a point $A_1$. Let us fix $\Delta r$ and among all rays $\gamma_{t,\Delta r}$,
$t\in[0,1]$, chose the one $\gamma_{t(\Delta r),\Delta r}$ that is passing through the point $A_1$.
Let us fix the point $O_1\in Oy$ by constructing the ball  $B(O_1,R+\Delta R(\Delta r_0))$ tangent
to the ray $\gamma_{1,\Delta r_0}$ and line $\xi$ passing through the point $A_1$.  (see Figure 3).
\begin{center}
\includegraphics[width=7 cm]{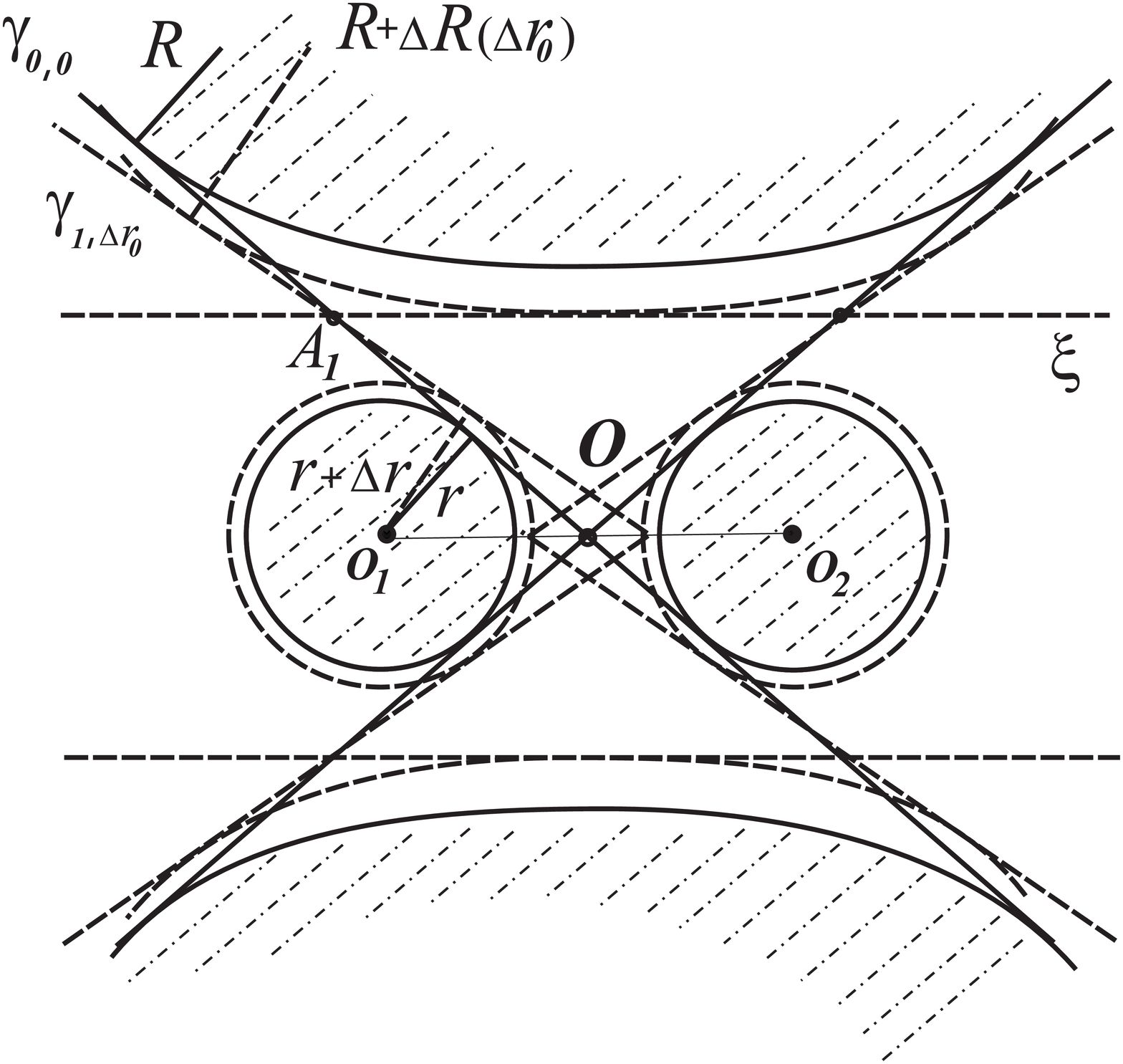}
\end{center}
\small\begin{center} Fig. 3 \end{center}
\normalsize
Now we can construct the system of concentric circles $\left\{ B(O_1,R+\Delta R(\Delta r))\right\}$ that are tangent to the respective rays $\gamma_{t(\Delta r),\Delta r}$ and the system $\left\{ B(O_2,R+\Delta R(\Delta r))\right\}$ symmetric to the first one with respect to the origin. It is easy to see that $\Delta R({\Delta r}_1)<\Delta R({\Delta r}_2)$ for any ${\Delta r}_1,{\Delta r}_2\in(0,{\Delta r}_0]$, such that ${\Delta r}_1<{\Delta r}_2$. Thus, the set of four closed balls $B:=\left\{ \overline{B(o_j,r)}, \overline{B(O_j,R)},\, j=1,2\right\}$ is approximated from the outside by the family of open sets $B_{\Delta r}:=\left\{ B(o_j,r+\Delta r), B(O_j,R+\Delta R(\Delta r)),\, j=1,2\right\}$, $0<\Delta r\le\Delta r_0$, and each $B_{\Delta r}$ is weakly $1$-semiconvex. With that, any ray starting at the origin $O\in \mathbb{R}^2\setminus B$ intersects set $B$. Thus, $B$ is weakly $1$-semiconvex and non $1$-semiconvex.
\end{example}

\begin{definition}\label{def5} A point  $x\in \mathbb{R}^n\setminus A$ is called {\it a
point of  $m$-nonsemiconvexity of a set $A\subset\mathbb{R}^n$} if there is no $m$-dimensional
half-plan that has $x$ on its boundary and does not intersect $A$.
\end{definition}

\begin{definition}

We say that {\it a set $A\subset \mathbb{R}^n$ is projected from a point $x\in\mathbb{R}^n$ on a set $B\subset \mathbb{R}^n$} if any ray, starting at point $x$ and intersecting $A$,  intersects $B$ as well.

\end{definition}


\begin{lemma}\label{lemm3} Let a set $E\subset \mathbb{R}^2$ be weakly $1$-semiconvex but non $1$-semiconvex and consist of three components. Then none of its components is projected on the union of the others from a point of $1$-nonsemiconvexity of $E$.
\end{lemma}

 \begin{theorem}\label{theor3}  Let  $E$ be open, bounded,  weakly $1$-semiconvex, and non $1$-semiconvex set with smooth boundary in $\mathbb{R}^2$. Then $E$ consists of minimum four components.
\end{theorem}

 \begin{lemma}\label{pro2}  Let a set $E\subset\mathbb{R}^n$ be closed, weakly $1$-semiconvex, and non $1$-semiconvex. Then for any family of open, weakly $1$-semiconvex sets $E_k$,
$k=1,2,\ldots$, approximating set $E$ from the outside, there exists $k_0\ge 1$ such that every set $E_k$,
$k=k_0,k_0+1,\ldots$, of the family is non $1$-semiconvex.
\end{lemma}

\begin{theorem}\label{theortheor} Let a set $E\subset \mathbb{R}^2$ be closed, weakly $1$-semiconvex and non $1$-semiconvex. Then $E$ consists of minimum three components.
\end{theorem}

\begin{theorem}\label{theor4}  Let  $E$ be closed, bounded,  weakly $1$-semiconvex, and non $1$-semiconvex set with smooth boundary in $\mathbb{R}^2$. Then $E$ consists of minimum  four components.
\end{theorem}

\begin{corollary} Let a set $E\subset \mathbb{R}^2$ be bounded, weakly $1$-semiconvex, non $1$-semiconvex, and consist of four components with smooth boundary. Then none of its components is projected on the union of the others from a point of $1$-nonsemiconvexity of $E$.
\end{corollary}

\smallskip

\noindent{\bf Acknowledgements.} {\it The author was supported by the grant of the President of
Ukraine for competitive projects F75//30308 of the State Fund for Fundamental Research.}

\end{document}